\documentclass[12pt]{amsart}

\usepackage{latexsym}
\usepackage{amssymb}
\usepackage{amsmath}
\usepackage{amscd}
\usepackage{amsthm}
\usepackage[all]{xy}  

\numberwithin{equation}{section}

\newtheoremstyle{fancy1}{10pt}{10pt}{\itshape}{12pt}{\textsc\bgroup}{.\egroup}{8pt}{
}
\newtheoremstyle{fancy2}{10pt}{10pt}{}{12pt}{\itshape}{.}{8pt}{ }

\theoremstyle{fancy1}

\newtheorem*{thethm*}{Theorem A}

\newtheorem{main}{Theorem}
\newtheorem*{main*}{Result}

\newtheorem*{cor*}{Corollary}

\theoremstyle{fancy2}

\newtheorem*{def*}{Definition}

\newtheorem*{rem*}{Remark}
\newtheorem*{rems*}{Remarks}

\newtheorem*{example*}{Example}

\newtheorem*{examples*}{Examples}
\newtheorem*{conclusion*}{Conclusion}

\theoremstyle{remark}

\newcommand{\cref}[1]{Corollary~\ref{#1}}

\newcommand{\RP}{\mathbb{R\mkern1mu P}}
\newcommand{\CP}{\mathbb{C\mkern1mu P}}

\newcommand{\bbone}{\mathbf{1}}

\newcommand{\R}{{\mathbb{R}}}
\newcommand{\Z}{{\mathbb{Z}}}
\renewcommand{\H}{{\mathbb{H}}}

\newcommand{\SO}{\ensuremath{\operatorname{\mathsf{SO}}}}

\newcommand{\Sp}{\ensuremath{\operatorname{\mathsf{Sp}}}}
\newcommand{\U}{\ensuremath{\operatorname{\mathsf{U}}}}
\newcommand{\SU}{\ensuremath{\operatorname{\mathsf{SU}}}}
\newcommand{\Spin}{\ensuremath{\operatorname{\mathsf{Spin}}}}

\renewcommand{\S}{\ensuremath{\operatorname{\mathsf{S}}}}

\def\con#1=#2(#3){#1 \equiv #2 \bmod{#3}}

\newcommand{\set}[1]{\left\{#1\right\}}
\newcommand{\grp}[1]{\langle #1 \rangle}
\newcommand{\ab}{\text{ab}}  

\newcommand{\diag}{\ensuremath{\operatorname{diag}}}

\renewcommand{\Im}{\ensuremath{\operatorname{Im}}}

\newcommand{\del}{\partial}

 \DeclareMathOperator{\pr}{pr}

\DeclareMathOperator{\Id}{Id} 
 
\DeclareMathOperator{\lcm}{lcm}
\DeclareMathOperator{\Span}{span}

\newcommand{\tmu}{\tilde{\mu}}
\newcommand{\tnu}{\tilde{\nu}}
\newcommand{\tu}{\tilde{u}}

\newcommand{\td}{\widetilde{d}}
\newcommand{\tA}{\widetilde{A}}
\newcommand{\tB}{\widetilde{B}}

\newcommand{\tD}{\widetilde{D}}
\newcommand{\tE}{\widetilde{E}}
\newcommand{\tH}{\widetilde{H}}
\newcommand{\tK}{\widetilde{K}}
\newcommand{\tL}{\widetilde{L}}

\newcommand{\bd}{\overline{d}}

\newcommand{\bA}{\overline{A}}
\newcommand{\bB}{\overline{B}}

\newcommand{\bD}{\overline{D}}
\newcommand{\bE}{\overline{E}}

\newcommand{\brho}{\bar{\rho}}

\newcommand{\hppm}{\widehat{p}_\pm}
\newcommand{\hqpm}{\widehat{q}_\pm}
\newcommand{\hpp}{\widehat{p}_+}
\newcommand{\hqp}{\widehat{q}_+}
\newcommand{\hpm}{\widehat{p}_-}
\newcommand{\hqm}{\widehat{q}_-}

\newcommand{\pom}{^+\!\!\!\!\diagup\!\!\!_-}
\newcommand{\mop}{^-\!\!\!\!\diagup\!\!\!_+}

\begin{document}


\title[Low Dimensional Cohomogeneity One Manifolds]{On the homology of low dimensional\\ cohomogeneity one manifolds}

\author{Corey A. Hoelscher}
\address{Rutgers University\\
New Brunswick, NJ}

\begin{abstract}
In this paper we give a characterization of the possible homology groups that can occur for compact simply connected cohomogeneity one manifolds in dimensions seven and lower.
\end{abstract}

\maketitle

%
%
%
%
%
%

\renewcommand{\thetable}{\Roman{table}}

\section*{Introduction}



One way to understand the size of the symmetry group of a manifold is by looking at the dimension of the orbits of the symmetry group. Heuristically, homogeneous spaces are the most symmetric in this sense, and cohomogeneity one manifolds are defined to be the next most symmetric. More precisely, a cohomogeneity one manifold is a connected smooth manifold with the smooth action of a compact connected Lie group with at least one orbit of codimension one, or equivalently with a one dimensional orbit space.

In addition to their natural importance to the theory of group actions on manifolds, cohomogeneity one manifolds are important in many areas of geometry. In the area of Riemannian geometry, it was shown in \cite{GZ2000} that a large class of cohomogeneity one manifolds admit metrics of non-negative sectional curvature. In particular, it was shown in \cite{thesis} that all but possibly two families of compact simply connected cohomogeneity one manifolds in dimensions 7 and lower, admit metrics of non-negative sectional curvature. Cohomogeneity one manifolds also give examples of non-negative and positive Ricci curvature \cite{GZ2002}. More recently, in \cite{GVZexample}, a cohomogeneity one structure was used to show that a manifold homeomorphic to the unit tangent bundle of $S^4$ admits a metric of positive sectional curvature. Cohomogeneity one manifolds are also important in mathematical physics as they give new examples of Einstein and Einstein-Sasaki manifolds (see \cite{Conti} or \cite{GHY}) and examples of manifolds with $G_2$ and $\Spin(7)$-holonomy (see \cite{CS} and \cite{CGLP}).

Cohomogeneity one manifolds were classified in dimensions 4 and lower in \cite{Ne} and \cite{Parker} and compact simply connected cohomogeneity one manifolds were classified in dimensions 5, 6 and 7 in \cite{thesis}. These classifications describe the possible manifolds by expressing them as the union of two disk bundles glued along their common boundary. However, based on this description, the topological invariants of the manifolds are not always clear. In lower dimensions, a description of precisely which manifolds occur up to diffeomorphism is known in the compact simply connected case. This was done for dimension 3 in \cite{Ne}, dimension 4 in \cite{Parker}, dimension 5 in \cite{thesis}, and recently dimension 6 in \cite{6dim_cohom1}. 


To give a similar characterization in dimension 7 would be much more difficult because the manifolds become more complicated and the required invariants are harder to compute. The first step in this direction is to compute the homology groups of these manifolds. From the classification \cite{thesis}, among compact simply connected cohomogeneity one manifolds, there are 4 primitive families and 9 non-primitive families in dimension 7 which have not been explicitly identified. The homology groups for the primitive families were computed in \cite{GWZ} and \cite{Ultman}. In this paper we compute the homology groups of the non-primitive manifolds, up to a group extension problem in a few cases.

\begin{main}\label{mainthrm}
If $M$ is a compact simply connected cohomogeneity one manifold of dimension 7 or less, then either $M$ has the same set of homology groups as a compact symmetric space $S$, or $M$ is 7-dimensional and has one of the exceptional sets of homology groups from Table \ref{t:homgrps}. Furthermore, in Table \ref{t:homgrps}, if $\alpha\ne 0$ then $\beta\in \set{1,\gamma}$.
\end{main}

It is known that all compact simply connected symmetric spaces in dimensions 7 or less admit cohomogeneity one actions \cite[Sec. 5.1]{thesis}. So their homology groups certainly occur as possibilities in Theorem \ref{mainthrm}.

Notice also that if $M$ is 7-dimensional and simply connected then all its homology groups are determined by $H_2(M)$ and $H_3(M)$ by Poincar\'{e} duality and the universal coefficients theorem.

In this paper we will use the notation $\Z/\delta=\Z/\delta\Z=\Z_\delta$ so for example $\Z/1=0$ and $\Z/0=\Z$. In particular, the last sentence in Theorem \ref{mainthrm} says that if $\alpha\ne 0$ then either $H_3(M)\simeq \Z/\gamma$ or we have a short exact sequence $0\to\Z/\gamma\to H_3(M)\to\Z/\gamma\to 0$.

{\setlength{\tabcolsep}{0.2cm}
\renewcommand{\arraystretch}{1.6}
\stepcounter{equation}
\begin{table}[!h]
\begin{center}
\begin{tabular}{|c||c|c|}
\hline
Type & $\qquad H_2(M) \qquad$ & \qquad $H_3(M)$ \qquad\qquad\\
\hline \hline
1 & $\quad H_2(M)=0\quad$ & $\quad H_3(M)=\Z/\gamma \quad$\\
\hline
2 & $\quad H_2(M)=\Z\oplus\Z/\alpha \quad$ &  $\quad 0\to\Z/\beta\to H_3(M)\to\Z/\gamma\to 0 \quad$ \emph{(exact)}\\
 &  where $\alpha\in\set{0,1,2}$& where $\beta\in\set{0,1,\dots}$ and $\gamma\in\set{1,2,\dots}$ \\
\hline
\end{tabular}
\end{center}
\vspace{0.1cm}
\caption{Exceptional sets of homology groups for $M^7$, from Theorem \ref{mainthrm}. If $\alpha\ne 0$ then $\beta\in \set{1,\gamma}$.}\label{t:homgrps}
\end{table}}


The paper is organized as follows. In the first section, we review the basic structure of cohomogeneity one manifolds and outline the basic tools that are used in computing the homology groups. The second section is the heart of the paper where we consider each remaining family of manifolds appearing in the classification \cite{thesis} and compute its homology groups. A subsection is devoted to each family separately. Some of these are easy to compute using the description of the manifolds as the union of two disk bundles, however others are quite difficult and require intricate techniques. The full descriptions of which manifolds have which homology groups are given in the conclusions at the end of Subsections \ref{N^7_A} through \ref{N^7_I}. There, we also give the precise values of the constants appearing in Table \ref{t:homgrps}. The reader who is interested in these results could skip the the ends of these sections.

The author would like to thank S.~K.~Ultman for several helpful discussions and for introducing him to the very useful results in \cite{Ultman}. He would also like to thank W.~Ziller for helpful discussions and for suggesting improvements to the manuscript.

%
%
%
%
%
%

\renewcommand{\thetable}{\theequation}

\section{Basic Techniques}\label{prelims}

In this section we discuss some tools and techniques that are used in the next section to compute the homology groups of the cohomogeneity one manifolds from the classification. First let us briefly recall the basic structure of a cohomogeneity one manifold. See \cite{thesis}, \cite{GZ2000} or \cite{GWZ} for more detailed descriptions. Suppose $G$ is a compact connected Lie group which acts by cohomogeneity one on a compact connected manifold $M$, with finite fundamental group. It then follows that the orbit space $M/G$ is an interval, say $[-1,1]$. It is clear that $M$ can be decomposed as the union of $D(B_-)=\pi^{-1}([-1,0])$ and $D(B_+)=\pi^{-1}([0,1])$, where $\pi:M\to M/G$ is the projection. Fixing a $G$-invariant Riemannian metric on $M$, the slice theorem says that $D(B_\pm)$ is a disk bundle over $B_\pm=\pi^{-1}(\pm1)$, with each point in $D(B_\pm)$ mapped to its closest point in $B_\pm$. Choose $x_0\in \pi^{-1}(0)$ and let $x_\pm$ be the closest point to $x_0$ in $B_\pm$. Then $D(B_\pm)$ is $G$-equivariantly diffeomorphic to $G\times_{K^\pm}D_\pm$ where $K^\pm$ is the isotropy group at $x_\pm$ and $D_\pm$ is the fiber at $x_\pm$ of the disk bundle $D(B_\pm)$, and where $K^\pm$ acts on $D_\pm$ via the slice representation. If we let $H=G_{x_0}$ be the isotropy subgroup of $G$ at $x_0$, then $H$ is also the isotropy group at $x_0$ for the $K^\pm$ action on $D_\pm$. In addition, $K^\pm$ acts transitively on $\del D_\pm$, so in fact, $K^\pm/H\simeq\del D_\pm=S^{\ell_\pm}$ is a sphere. In conclusion we can describe $M$, $G$-equivariantly, as
\begin{equation}\label{decomp}
  M\,\simeq\, G\times_{K^-}D_- \,\, \cup \,\, G\times_{K^+}D_+ \text{ \hspace{1em} where \hspace{1em} } S^{\ell_\pm}=\del D_\pm \simeq K^\pm/H
\end{equation}
and where the two halves $G\times_{K^\pm}D_\pm$ are glued along their common boundary $G\times_{K^\pm}K^\pm/H=G/H$. Hence $M$ is described entirely in terms of the isotropy groups $K^\pm$ and $H$. Conversely, given compact groups $G\supset K^-,K^+\supset H$ with $K^\pm/H\simeq S^{\ell_\pm}$. We can build a cohomogeneity one manifold $M$ using \eqref{decomp}. This collection of groups $G\supset K^-,K^+\supset H$ is called the \emph{group diagram} of $M$.



\subsection{The Mayer-Vietoris sequence}

Given the decomposition of the cohomogeneity one manifold $M=G\times_{K^-}D_- \cup G\times_{K^+}D_+$ described above, the simplest way to attempt to compute the cohomology groups of $M$ is through the Mayer-Vietoris sequence. Notice that $G\times_{K^\pm}D_\pm$ deformation retracts to $G/K^\pm$ and this retraction takes the boundary $G/H$ to $G/K^\pm$ via the standard projection $\pi_\pm:G/H\to G/K^\pm$. Hence we have the long exact sequence
\begin{equation}\label{MVseq}
\begin{CD}
\cdots\to H^n(M) @>{(i_-^*,i_+^*)}>> H^n(G/K^-)\oplus H^n(G/K^+) @>{\pi_-^* - \pi_+^*}>>H^n(G/H) \to H^{n+1}(M)\to\cdots
\end{CD}
\end{equation}
where $i_\pm:G/K^\pm\to M$ is the inclusion.

\subsection{The long exact sequence of the pair}

Recall $K^-/H\simeq S^{\ell_-}$ and suppose the bundle $K^-/H \to G/H \to G/K^-$ is orientable as a sphere bundle \cite[pg. 442]{Hatcher}. In this case, the author of \cite[Sec. 4]{Hebda} describes how the long exact sequence for the pair $(M,B_+)$ can be modified using the Thom isomorphism to give the following long exact sequence:
\begin{equation}\label{lespair}
\begin{CD}
\cdots\to H^{\,n-\ell_- -1}(G/K^-) \to H^n(M) @>{i_+^*}>>H^n(G/K^+) \to H^{\,n-\ell_-}(G/K^-)\to\cdots
\end{CD}
\end{equation}
and similarly for $(M,B_-)$ if $S^{\ell_+}\to G/H \to G/K^+$ is an orientable sphere bundle. This sequence was first used for cohomogeneity one manifolds in \cite{Ultman}.

\subsection{Non-primitive actions}\label{sec:nonprim}


Recall that a cohomogeneity one manifold $M$ is called non-primitive if the action has a group diagram $G\supset K^-,K^+\supset H$ such that there is a compact connected proper subgroup $L\subset G$ which contains $K^-$, $K^+$ and $H$. If $M_L$ is the cohomogeneity one manifold given by the group diagram $L\supset K^-, K^+\supset H$, then $M$ is $G$-equivariantly diffeomorphic to $G\times_L M_L$, and we have the fiber bundle
\begin{equation}\label{nonprim}
M_L \to M \to G/L.
\end{equation}
See \cite{thesis} for more details. We will refer to this bundle as the \emph{non-primitivity fiber bundle}. In particular $L$ is the structure group for this bundle, though not necessarily effectively.

%
%
%
%
%
%

\section{Computing the Homology Groups}\label{body}

In this section we will prove the main theorem from the introduction. Suppose that $M$ is a compact simply connected cohomogeneity one manifold of dimension 7 or less. We will show that $M$ has the homology groups of a symmetric space or one of the exceptional sets of homology groups listed in Table \ref{t:homgrps}. 

First suppose $\dim(M)\le 4$. The only compact simply connected cohomogeneity one manifolds here are $S^2$, $S^3$, $S^4$, $S^2\times S^2$, $\CP^2$ and $\CP^2\# -\CP^2$ (see \cite{Parker}). It is clear these all have the homology groups of a symmetric space. If $\dim(M)=5$, we know from \cite[Thm. C]{thesis} that $M$ is diffeomorphic to $S^5$, $\SU(3)/\SO(3)$, $S^3\times S^2$ or the nontrivial $S^3$ bundle over $S^2$. It is also clear that these manifolds have the homology groups of a symmetric space.

In the case that $M$ is 6-dimensional, we know from \cite{6dim_cohom1} that $M$ is diffeomorphic to a symmetric space; an $S^2$ bundle over $\CP^2$ or $S^2\times S^2$; an $S^4$ bundle over $S^2$; or a manifold of type $N^6_D$. The first case is a tautology. In the next two cases, the Gysin sequence clearly shows that the homology groups of $M$ are the same as those of $S^2\times \CP^2$, $S^2\times S^2\times S^2$ or $S^2\times S^4$. In the last case, $N^6_D$, $M$ is given by the group diagram $S^3\times S^3\supset T^2, S^3\times S^1\supset \set{(e^{ip\theta}, e^{i\theta})}$. Here $G/K^-\simeq S^2\times S^2$, $G/K^+\simeq S^2$, and $G/H=S^3\times S^3/\set{(e^{ip\theta}, e^{i\theta})}\simeq S^3\times S^2$ (see \cite[Prop. 2.3]{WZ}). The Mayer-Vietoris sequence \eqref{MVseq} easily shows $H^4(M)\simeq \Z^2$ and therefore $H^3(M)$ must be torsion free and $H^2(M)\simeq \Z^2$. Then there is a segment of \eqref{MVseq} of the form $\Z^2\to \Z^3\to \Z\to H^3(M)\to 0$ and hence $H^3(M)=0$. Therefore $M$ has the same homology groups as $\CP^2\times S^2$ in this case.

Finally suppose $\dim(M)=7$. We know from \cite[Thm. A]{thesis} that $M$ is diffeomorphic to one of the following: a symmetric space; a Brieskorn variety $B_d^7$; the product of a lower dimensional cohomogeneity one manifold with a homogeneous space; or a manifold given by one of the group diagrams listed in Tables I and II of \cite{thesis}. The first case is obvious. Next, if $M\simeq B_d^7$ then the homology groups of $M$ are know to be given by $H_2(M)=0$ and $H_3(M)\simeq \Z/d$ \cite[Cor. V9.3]{Br}, and these groups appear in Table \ref{t:homgrps}.

Next suppose $M=N_1\times N_2$ is the product of a cohomogeneity one manifold with a homogeneous space. Either $N_1$ or $N_2$ must be a sphere, since spheres are the only compact simply connected homogeneous or cohomogeneity one manifolds in dimensions 3 or lower. Say $M=S^k\times N$, where $N$ is either homogeneous or cohomogeneity one, with dimension 5 or less. It is clear from above, and from the classification of low dimensional homogeneous spaces (e.g.~as in \cite[Prop. 2.1]{thesis}), that $N$ must have the homology groups of a symmetric space. By the K\"{u}nneth formula, $M$ must also have the homology of a symmetric space. 

We are only left with the case that $M$ is given by one of the group diagrams from Tables I and II of \cite{thesis}. The homology of the manifolds of type $P^7_A$ and $P^7_D$ was computed in \cite{Ultman}. They showed, for $P^7_A$, that $H_3(M)\simeq \Z/r$ and either $H_2(M)\simeq \Z$ or $H_2(M)\simeq \Z\oplus \Z_2$, and that $H_3(M)$ is finite (i.e.~$r>0$) in the case $H_2(M)\simeq \Z\oplus \Z_2$. Also for $P^7_D$, they showed $H_2(M)\simeq \Z$ and $H_3(M)\simeq \Z/r$. Next, the manifolds of type $P^7_B$ and $P^7_C$, were considered in \cite[Sec. 13]{GWZ}. They showed $H_2(M)\simeq \Z$ and $H_3(M)\simeq \Z/r$ for $M$ of type $P^7_B$; and $H_2(M)=0$ and $H_3(M)\simeq \Z/r$ for $M$ of type $P^7_C$. In each case, when $r>0$, we see these sets of homology groups appear in Table \ref{t:homgrps} above. If $r=0$ then these homology groups match those of a symmetric space.

The only manifolds which are left are those 7 dimensional manifolds appearing in Table II of \cite{thesis}. For the rest of this section we consider these manifolds, one by one, and compute their homology groups. In each case we start by recalling the group diagram for the action and the conditions on the diagram which make the resulting manifold $M$ simply connected. We conclude each case with a review of the homology groups of $M$.

The group diagrams are given in the form $G\supset K^-, K^+\supset H$ and $H=H_-\cdot H_+$ in each case where $H_\pm=H\cap K^\pm_0$. Note also that $i,j,k\in S^3$ denote the standard unit quaternions and when we write $\set{e^{i\theta}}\subset S^3$ we understand this to mean $\set{e^{i\theta}=\cos\theta+i\sin\theta \,|\, \theta\in\R}$ and similarly $\set{z}\subset S^3$ will mean $\set{z \,|\, z=e^{i\theta}\in S^3, \theta\in \R}$.

\subsection{Actions of type $N^7_A$:}\label{N^7_A}

$$S^3\times S^3 \,\, \supset \,\, \set{(e^{ip_-\theta},e^{iq_-\theta})}\cdot H_+, \,\, \set{(e^{ip_+\theta},e^{iq_+\theta})} \cdot H_- \,\, \supset \,\,  H_-\cdot H_+$$
where $H_\pm=H\cap K^\pm_0$ is finite cyclic and $\gcd(p_\pm,q_\pm)=1$.

\vspace{.5em}

We divide this family into two cases, depending on whether or not $K^-_0=K^+_0$. First suppose $K^-_0=K^+_0$. Since $H=H_-\cdot H_+$ it follows that $H\subset K^\pm_0$ and hence $K^\pm=K^\pm_0$. Then we know $G/K^\pm\simeq S^3\times S^2$ by \cite[Prop. 2.3]{WZ}. Since $G/K^\pm$ is simply connected, $K^-/H\to G/H\to G/K^-$ is an orientable sphere bundle. Hence we have the long exact sequence of the pair \eqref{lespair}. This sequence clearly give $H^5(M)\simeq \Z\oplus \Z$. Next consider the non-primitivity fiber bundle \eqref{nonprim} with $L=K^\pm$. This takes the form $S^2\to M\to S^3\times S^2$. The Gysin sequence for this bundle clearly gives $H^4(M)\simeq \Z$, using the fact that $H^5(M)\simeq \Z\oplus \Z$. So if $K^-_0=K^+_0$ then $M$ has the homology groups of $S^3\times S^2\times S^2$. In fact it follows that the Euler class of this bundle must be trivial.

The case where $K^-_0\ne K^+_0$ will be much more complicated. We start by establishing some notation. Denote $b_\pm:= |H_\pm|=|H\cap K^\pm_0|$ and $h:=|H|$ so that $K^\pm$ has $h/b_\pm$ connected components. Further let $a=q_+p_--q_-p_+$ and notice that the number of intersection points of $K^-_0$ and $K^+_0$ is $|a|$.

Let us start by looking at the non-primitivity bundle \eqref{nonprim} with $L=T^2$. We see that $G/L\simeq S^2\times S^2$, and that the effective version of the group diagram for $M_L$ is $L/H \supset K^-/H, K^+/H \supset 1$. We then identify this as the group diagram for a lens space $S^3/\Z_r$, where $r$ is the order of the intersection of the two circles $K^-/H$ and $K^+/H$, \cite[Sec. 7.2]{thesis}. It is not too difficult to see that $r=|ah/b_-b_+|$.

Now consider the spectral sequence for the non-primitivity bundle $S^3/\Z_r \to M\to S^2\times S^2$. The first page takes the form $E^{p,q}_2= H^p(S^2\times S^2)\otimes H^q(S^3/\Z_r)$ with maps $d_2^{p,q}:E^{p,q}_2\to E_2^{p+2,q-1}$. First notice that $E_\infty^{2,3} = E_3^{2,3}= \ker(d_2^{2,3}: \Z\oplus\Z\to \Z_r)\simeq \Z\oplus\Z$. Hence $H^5(M)\simeq \Z\oplus\Z$, since $E_\infty^{2,3}$ is the only nontrivial group along that diagonal.

Next, we claim that $E_\infty^{2,2}\simeq \Z/r$. For this, first recall $H^6(M)=0$ since $M^7$ is simply connected. In particular $E_\infty^{4,2}$ must be trivial which means $d_2^{2,3}:E^{2,3}_2\to E_2^{4,2}$ must be onto. It follows that there are elements $w\in E_2^{0,3}$ and $v\in E_2^{2,0}$ such that $d_2^{2,3}(vw)$ generates $E_2^{4,2}\simeq \Z/r$. Then $d_2^{2,3}(vw)=d_2^{2,0}(v)w+vd_2^{0,3}(w)=vd_2^{0,3}(w)$. Since this element generates $E_2^{4,2}\simeq \Z/r$ it follows that $d_2^{0,3}(w)\in E_2^{2,2}\simeq \Z_r\oplus \Z_r$ must have order $r$. Hence the image of $d_2^{0,3}:\Z\to \Z_r\oplus \Z_r$ must be cyclic of order $r$. It follows that $E_\infty^{2,2}= E_3^{2,2}= E_2^{2,2}/\Im(d_2^{0,3}) \simeq \Z_r$. The only other nontrivial group on the same diagonal as $E_\infty^{2,2}$ is $E_\infty^{4,0}= E_4^{4,0}/\Im(d_4^{0,3}: \Z \to \Z)\simeq \Z/\beta$ for some $\beta\in \set{0,1,2,3,\dots}$. Therefore $H^4(M)$ fits into the short exact sequence $0\to \Z/\beta \to H^4(M)\to \Z/\gamma$, for $\gamma=r=|ah/b_-b_+|$. The task for the rest of this section is to compute the value of $\beta$.

To find $\beta$ we only need the order of $H^4(M)$ and we will eventually compute this using the Mayer-Vietoris sequence \eqref{MVseq}. Before we can do this we need to compute $H^i(G/K^\pm)$, $H^i(G/H)$ and the maps $H^i(G/K^\pm)\to H^i(G/H)$ induced from the projections $\pi_\pm: G/H\to G/K^\pm$. We will accomplish these tasks by considering the following diagram of fiber bundles:
\begin{equation}\label{pidiagram}
\displaystyle{\begin{array}{ccccc}
      L/H & \to & G/H & \to & G/L\\
      \downarrow & & \downarrow & & \downarrow \\
      L/K^\pm & \to & G/K^\pm & \to & G/L
 \end{array}  }
\end{equation}
where $L=T^2\subset G$, as before. The middle vertical map is $\pi_\pm$. $G/L\simeq S^2\times S^2$ and the vertical map $G/L\to G/L$ is the identity.

The top bundle gives a spectral sequence with $E_2^{p,q}=H^p(G/L)\otimes H^q(L/H)$ and maps $d_2^{p,q}:E_2^{p,q}\to E_2^{p+2,q-1}$. The bottom bundle also has a spectral sequence $\bE_2^{p,q}=H^p(G/L)\otimes H^q(L/K^\pm)$ and maps $\bd_2^{p,q}:\bE_2^{p,q}\to \bE_2^{p+2,q-1}$. We will soon compute all these groups and maps explicitly, as well as the maps from $\bE_2$ to $E_2$ induced by $\pi_\pm$. However, we must first find a concrete set of fixed generators of $E_2^{0,1}\simeq H^1(L/H)\simeq \Z^2$.

To understand $L/H$ choose a homomorphism $\rho: L\to T^2$ with $\ker(\rho)=H$ and suppose $\rho(z,w)=(z^Aw^B,z^Dw^E)$. Then $\rho(K^\pm)= \rho(K^\pm_0)= \rho(\set{(z^{p_\pm},z^{q_\pm})})= \set{(z^{Ap_\pm+Bq_\pm},z^{Dp_\pm+Eq_\pm})}$. Furthermore, since $\ker(\rho)\cap K^\pm_0=H\cap K^\pm_0\simeq \Z_{b_\pm}$, we see $\gcd(Ap_\pm+Bq_\pm,Dp_\pm+Eq_\pm)=b_\pm$, so we can define the integers $\hppm:=(Ap_\pm+Bq_\pm)/b_\pm$ and $\hqpm:=(Dp_\pm+Eq_\pm)/b_\pm$. With this notation $ \rho(K^\pm)= \set{(z^{\hppm},z^{\hqpm})}$ and $\gcd(\hppm,\hqpm)=1$. Also notice that since $\ker(\rho)=H$, $|AE-BD|=|H|=h$.

For example, one choice of $\rho$ can be constructed as follows. First we establish some notation which will simplify our description. For every fixed ordered pair $(x,y)\in \Z^2$, make a choice of relatively prime integers $\phi(x,y),\psi(x,y)$ such that $x\psi_{(x,y)}-y\phi_{(x,y)}=\gcd(x,y)$. Also consider the standard projection $\R^2\to T^2:(s,t)\mapsto (e^{2\pi is},e^{2\pi it})$. Every homomorphism $f:T^2\to T^2$ has a unique lift to a linear map $\bar f:\R^2\to \R^2$, via this cover. With this notation, first define $\displaystyle{\brho_1= \left[\begin{array}{cc} \psi(p_-,q_-) & -\phi(p_-,q_-) \\ -q_- & p_-  \end{array}\right] }$ and $\displaystyle{\brho_2= \left[\begin{array}{cc} b_- & 0 \\ 0 & 1  \end{array}\right] }$. Then notice $\rho_1$ is an isomorphism which maps $K^-_0$ to $\set{(z,1)}$. Hence $\ker(\rho_2\circ\rho_1)=(\Z/b_-\subset K^-_0)=K^-_0\cap H=H_-$. Furthermore $\rho_2\circ\rho_1(H)\simeq H/H_-\simeq H_+/(H_-\cap H_+)$ so the order of this set is $h/b_-=b_+/|H_-\cap H_+|$. In particular $|H_-\cap K^+_0|=|H_-\cap H_+|=b_-b_+/h$ and so $\ker(\rho_2\circ\rho_1|_{K^+_0})$ is cyclic of order $b_-b_+/h$. We also see, for $(z^{p_+},z^{q_+})\in K^+_0$, that $\rho_2\circ\rho_1(z^{p_+},z^{q_+}) = (z^{b_-c},z^a)$ where $c=p_+\psi(p_-,q_-)-q_+\phi(p_-,q_-)$. It then follows that $\gcd(b_-c,a)=b_-b_+/h$. Then we can write $\rho_2\circ\rho_1(K^+_0)=\set{(z^{hc/b_+},z^{ha/b_-b_+})}$ with $\gcd(hc/b_+,ha/b_-b_+)=1$. Now define $\displaystyle{\brho_3= \left[\begin{array}{cc} \psi(\frac{hc}{b_+},\frac{ha}{b_-b_+}) & -\phi(\frac{hc}{b_+},\frac{ha}{b_-b_+}) \\ -\frac{ha}{b_-b_+} & \frac{hc}{b_+}  \end{array}\right] }$ and $\displaystyle{\brho_4= \left[\begin{array}{cc} h/b_- & 0 \\ 0 & 1  \end{array}\right] }$. As before $\rho_3$ is an isomorphism which maps $\rho_2\circ\rho_1(K^+_0)$ to $\set{(z,1)}$. Also, the kernel of $\rho_4\circ\rho_3$ is $\Z/(h/b_-)\subset \rho_2\circ\rho_1(K^+_0)$. Yet we also know that $\rho_2\circ\rho_1(H) =\Z/(h/b_-)\subset \rho_2\circ\rho_1(K^+_0)$. Putting this together with the facts that $\ker(\rho_2\circ\rho_1)=H_-$ and $H=H_-\cdot H_+$ we conclude $\ker(\rho_4\circ\rho_3\circ\rho_2\circ\rho_1)=H$. Hence $\rho_4\circ\rho_3\circ\rho_2\circ\rho_1$ could be taken as our function $\rho$. With this $\rho$ we would then have
\begin{equation}\label{rhoeg}
\displaystyle{ \left[\begin{array}{cc} A & B \\ D & E  \end{array}\right] = \brho=\brho_4\brho_3\brho_2\brho_1= }
\end{equation}
$$
\displaystyle{\left[\begin{array}{cc} \frac{h}{b_-}\left(\psi_{(\frac{hc}{b_+},\frac{ha}{b_-b_+})}\psi_{(p_-,q_-)}b_- + \phi_{(\frac{hc}{b_+},\frac{ha}{b_-b_+})}q_-\right) & -\frac{h}{b_-}\left(\psi_{(\frac{hc}{b_+},\frac{ha}{b_-b_+})}\phi_{(p_-,q_-)}b_- + \phi_{(\frac{hc}{b_+},\frac{ha}{b_-b_+})}p_-\right) \\ -\frac{hq_+}{b_+} & \frac{hp_+}{b_+}  \end{array}\right]}.
$$
Furthermore, we would then have $\hpp=1$, $\hqp=0$, $\hpm=h\psi(\frac{hc}{b_+},\frac{ha}{b_-b_+})/b_-$ and $\hqm=-ha/b_-b_+$. However, because these are so complicated, we will continue to leave everything in terms of $A,B,D,E$.

We are now ready to choose our generators for $H^1(L/H)$ using $\rho$. Notice $\rho$ induces a natural isomorphism $f:L/H\to T^2$ given by $f(\ell H)=\rho(\ell)$. Now $H^1(T^2)$ is generated by $\pr_i^*(\bbone)$, for $i=1,2$, where $\pr_i:T^2\to S^1:(z_1,z_2)\mapsto z_i$ is the projection and where $\bbone$ is some fixed generator of $H^1(S^1)$. Define $u_1=f^*(\pr_1^*(\bbone))$ and $u_2=f^*(\pr_2^*(\bbone))$, as our fixed generators for $H^1(L/H)$.

We can now understand the map $\pi_\pm^*:H^1(L/K^\pm)\to H^1(L/H)$ in terms of these generators. Notice the projection $\tilde{\pi}_\pm: T^2 \to T^2/f(K^\pm/H)$ and the map $f:L/H\to T^2$ induce an isomorphism $\tilde{f}:L/K^\pm\to T^2/f(K^\pm/H)$, by $\tilde{f}(\pi_\pm(\ell H))=\tilde{\pi}_\pm(f(\ell H))$. Further $f(K^\pm/H)=\rho(K^\pm)=\set{(z^{\hppm},z^{\hqpm})}$ so there is a generator $\tilde{u}_\pm$ of $H^1(T^2/f(K^\pm/H))\simeq \Z$ such that $\tilde{\pi}_\pm^*(\tilde{u}_\pm)= \hqpm\pr_1^*(\bbone)- \hppm\pr_2^*(\bbone)$. Define $u_\pm=\tilde{f}^*(\tilde{u}_\pm)$ as our fixed generator of $H^1(L/K^\pm)\simeq \Z$. Then $\pi_\pm^*(u_\pm)= \hqpm u_1-\hppm u_2$.

We will now compute our map $\bd_2^{0,1}:\bE_2^{0,1}\to \bE_2^{2,0}$ which will then allow us to compute $H^i(G/K^\pm)$. For this consider the diagram of fiber bundles
\begin{equation}\label{hopfdiagram}
\displaystyle{\begin{array}{ccccc}
      L & \to & G & \to & G/L\\
      \downarrow & & \downarrow & & \downarrow \\
      L/F & \to & G/F & \to & G/L
 \end{array}  }
\end{equation}
with $F=K^\pm$. The top bundle is a product of Hopf bundles and has spectral sequence $\tE_2^{p,q}$ with maps $\td_2^{p,q}$. Let $\tilde{u}_1,\tilde{u}_2\in H^1(L)$ and $v_1,v_2\in H^2(G/L)$ be the standard product generators $\tilde{u}_i=\pr_i^*(\bbone)$ and $v_i=\pr_i^*(\bbone')$ for generators $\bbone\in H^1(S^1)$ as above and $\bbone'\in H^2(S^2)$. Then it is clear that with the right choice of sign on $\bbone'$ we have $\td_2^{0,1}(\tilde{u}_i)=v_i$. Now consider the projection map $L\to L/K^\pm$. Since $K^\pm_0=\set{(z^{p_\pm},z^{q_\pm})}$, we know the projection $L\to L/K^\pm_0$ induces a map $H^1(L/K^\pm_0)\to H^1(L)$ which sends a generator to $\pom(q_\pm \tu_1-p_\pm \tu_2)$. Then, since $L/K^\pm_0\to L/K^\pm$ is an $(h/b_\pm)$--fold cover of the circle, it follows the map $H^1(L/K^\pm)\to H^1(L)$ takes a generator to $\pom\frac{h}{b_\pm}(q_\pm \tu_1-p_\pm \tu_2)$. Then using the commutativity of the diagram of spectral sequences $\bE$ and $\tE$ from \eqref{hopfdiagram}, we see that $\bd_2^{0,1}(u_\pm) =\pom\frac{h}{b_\pm}(q_\pm v_1-p_\pm v_2)$.

Now we can compute the cohomology groups of $G/K^\pm$ from the spectral sequence $\bE$ described above. It is clear that $\bd_2^{0,1}$ is injective so we see $\bE_3^{0,1}=\bE_\infty^{0,1}=0$. Hence $H^1(G/K^\pm)=0$. Next, the image of $\bd_2^{0,1}$ in $\bE_2^{2,0}=\grp{v_1,v_2}_\text{ab}$ is generated by $\frac{h}{b_\pm}(q_\pm v_1-p_\pm v_2)$ so $\bE_\infty^{2,0}= \bE_3^{2,0}= \grp{v_1,v_2}_\text{ab}/\grp{\frac{h}{b_\pm}(q_\pm v_1-p_\pm v_2)} \simeq \Z\oplus(\Z/(h/b_\pm))$, since $p_\pm$ and $q_\pm$ are coprime. It follows $H^2(G/K^\pm)\simeq \Z\oplus(\Z/(h/b_\pm))$. Next, $\bE_2^{2,1}\simeq \Z\oplus \Z$ is generated by $v_1u_\pm$ and $v_2u_\pm$ and since $\bd_2^{2,0}(v_i)=0$ we have $\bd_2^{2,1}(v_iu_\pm) = v_i\bd_2^{0,1}(u_\pm) = \pom v_i\frac{h}{b_\pm}(q_\pm v_1-p_\pm v_2)$. In particular $\bd_2^{2,1}(v_1u_\pm) = \mop \frac{hp_\pm}{b_\pm}v_1 v_2$ and $\bd_2^{2,1}(v_2u_\pm) = \pom \frac{hq_\pm}{b_\pm}v_1 v_2$ so $\bE_\infty^{2,1}=\bE_3^{2,1}=\ker(\bd_2^{2,1})=\grp{q_\pm v_1u_\pm + p_\pm v_2u_\pm}$. This fact will be important later, but for now we simply conclude $H^3(G/K^\pm)\simeq \Z$. We also see $\Im(\bd_2^{2,1})=\grp{\frac{hp_\pm}{b_\pm}v_1 v_2, \frac{hq_\pm}{b_\pm}v_1 v_2}=\grp{\frac{h}{b_\pm}v_1 v_2}$ since $p_\pm$ and $q_\pm$ are coprime. So $\bE_\infty^{4,0}=\bE_3^{4,0}=\grp{v_1v_2}/\grp{\frac{h}{b_\pm}v_1 v_2}\simeq \Z/(h/b_\pm)$ and hence $H^4(G/K^\pm)\simeq \Z/(h/b_\pm)$. We also clearly see $H^5(G/K^\pm)\simeq \Z$ from $\bE$.

We now repeat this process for $G/H$. Again we have the diagram \eqref{hopfdiagram} with $F=H$ in this case. To understand the map induced by the projection $\wp: L\to L/H$ recall $f(\wp(\ell))=\rho(\ell)$ and recall $H^1(L/H)$ is generated by $u_i=f^*(\pr_i^*(\bbone))$, $i=1,2$. Therefore $\wp^*(u_i)= \rho^*(\pr_i^*(\bbone))$. Since $\rho(z,w)= (z^Aw^B,z^Dw^E)$ we see $\wp^*(u_1)= \rho^*(\pr_1^*(\bbone))= A\pr_1^*(\bbone)+B\pr_2^*(\bbone)= A\tu_1+B\tu_2$ and $\wp^*(u_2)= \rho^*(\pr_2^*(\bbone))= D\pr_1^*(\bbone)+E\pr_2^*(\bbone)= D\tu_1+E\tu_2$. It then follows from the commutativity of the diagram of spectral sequences from $E_2$ to $\tE_2$, that $d_2^{0,1}(u_1)=Av_1+Bv_2$ and $d_2^{0,1}(u_2)=Dv_1+Ev_2$.

We now compute $H^i(G/H)$ from the spectral sequence $E_s^{p,q}$. First observe that $E_s^{p,q}$ clearly shows $H^6(G/H)\simeq \Z$ so $G/H$ is orientable. Next notice that since $\rho(z,w)= (z^Aw^B,z^Dw^E)$ is onto, the vectors $Av_1+Bv_2$ and $Dv_1+Ev_2$ are linearly independent. Hence $d_2^{0,1}$ is injective, $E_\infty^{0,1}= E_3^{0,1}=0$, and hence $H^1(G/H)=0$. Next, $\Im(d_2^{0,1})= \grp{Av_1+Bv_2,Dv_1+Ev_2}$ so $E_\infty^{2,0}= E_3^{2,0}= \grp{v_1,v_2}_\ab / \grp{Av_1+Bv_2,Dv_1+Ev_2}$. Further $d_2^{0,2}(u_1u_2)= (Av_1+Bv_2)u_2-u_1(Dv_1+Ev_2) = -Dv_1u_1+Av_1u_2-Ev_2u_1+Bv_2u_2$, since $u_iv_j= (-1)^{1\cdot2}v_ju_i= v_ju_i$. In particular $d_2^{0,2}$ is injective so $E_\infty^{0,2}=E_3^{0,2}=0$. It then follows that $H^2(G/H) \simeq E_\infty^{2,0}= \grp{v_1,v_2}_\ab / \grp{Av_1+Bv_2,Dv_1+Ev_2}$. We then see $H^2(G/H) \simeq H$ as groups, either directly from this or from the obvious fact $\pi_1(G/H)\simeq H$. Next notice $d_2^{2,1}(v_iu_j)=v_id_2^{0,1}(u_j)$, since $d_2^{2,0}(v_i)=0$. So we have $d_2^{2,1}(v_1u_1)=Bv_1v_2$, $d_2^{2,1}(v_1u_2)=Ev_1v_2$, $d_2^{2,1}(v_2u_1)=Av_1v_2$ and $d_2^{2,1}(v_2u_2)=Dv_1v_2$. In particular $\Im(d_2^{2,1})= \grp{Av_1v_2,Bv_1v_2,Dv_1v_2,Ev_1v_2}= \grp{\delta v_1v_2}$ where $\delta=\gcd(A,B,D,E)$. Therefore $E_\infty^{4,0}= E_3^{4,0}= \grp{v_1v_2}/\grp{\delta v_1v_2} \simeq \Z/\delta$. We also see $E_\infty^{2,2}=\ker(d_2^{2,2}: \Z^2 \to \Z^2)$ is torsion free and hence must be trivial, since $H^4(G/H)$ is finite by Poincar\'{e} duality. It then follows that $H^4(G/H)\simeq E_\infty^{4,0}\simeq \Z/\delta$. By Poincar\'{e} duality again, the torsion subgroup of $H^3(G/H)$ is $\Z/\delta$. It then follows $H^3(G/H)\simeq E_\infty^{2,1}= E_3^{2,1}= \ker(d_2^{2,1})/\Im(d_2^{0,2})\simeq \Z^3/\Z\simeq \Z\oplus\Z\oplus(\Z/\delta)$.

We now have the cohomology groups for $G/K^\pm$ and $G/H$. Next we will determine the images of the maps $\pi_\pm^*: H^i(G/K^\pm)\to H^i(G/H)$, for $i=2,3,4$, for use in the Mayer-Vietoris sequence for $M$. Recall the diagram of bundles \eqref{pidiagram} which induces maps $(\pi_\pm)_s^{p,q}:\bE_s^{p,q}\to E_s^{p,q}$ on the spectral sequences. Since the map $G/L\to G/L$ in \eqref{pidiagram} is the identity, the maps $(\pi_\pm)_2^{i,0}$ are the identity as well. In particular $(\pi_\pm)_2^{2,0}(v_j)=v_j$. Then, since $\bE_\infty^{2,0}$ and $E_\infty^{2,0}$ are both quotients of $\grp{v_1,v_2}_\ab$ we see $(\pi_\pm)_\infty^{2,0}: [v_i]\mapsto [v_i]$ is obviously onto. Furthermore, since $H^2(G/K^\pm)\simeq \bE_\infty^{2,0}$ and $H^2(G/H)\simeq E_\infty^{2,0}$ in the natural way we see the map $H^2(G/K^\pm)\to H^2(G/H)$ is onto. In the same way $H^4(G/K^\pm)\to H^4(G/H)$ is onto. Next recall $\bE_\infty^{2,1}= \grp{q_\pm v_1u_\pm + p_\pm v_2u_\pm}$. We then see $(\pi_\pm)_\infty^{2,1}(q_\pm v_1u_\pm + p_\pm v_2u_\pm)= (q_\pm v_1 + p_\pm v_2)(\pi_\pm)_\infty^{0,1}(u_\pm)= (q_\pm v_1 + p_\pm v_2)(\hqpm u_1-\hppm u_2)=
q_\pm \hqpm v_1u_1 - q_\pm \hppm v_1u_2 + p_\pm \hqpm v_2u_1 - p_\pm \hppm v_2u_2$. Again since $H^3(G/K^\pm)\simeq \bE_\infty^{2,1}$ and $H^3(G/H)\simeq E_\infty^{2,1}$ in the natural way, we see that the image of $\pi_\pm^*$ in $H^3(G/H)$ is generated by the element $q_\pm \hqpm v_1u_1 - q_\pm \hppm v_1u_2 + p_\pm \hqpm v_2u_1 - p_\pm \hppm v_2u_2$.

Finally we have the Mayer-Vietoris sequence \eqref{MVseq} for $M$. The important segment of the sequence takes the form
\begin{equation}\label{ourMVseq}
\begin{CD}
@>{\Delta_2}>> H^3(M) @>{I_3}>> H^3(G/K^-)\oplus H^3(G/K^+) @>{\Pi_3}>> H^3(G/H) \\
@>{\Delta_3}>> H^4(M) @>{I_4}>> H^4(G/K^-)\oplus H^4(G/K^+) @>{\Pi_4}>> H^4(G/H) @>{\Delta_4}>>
\end{CD}
\end{equation}
where $\Pi_j=\pi_-^*-\pi_+^*$ and $I_j=(i_-^*,i_+^*)$. First notice that $\Delta_2=0$ and $\Delta_4=0$ since $\pi_\pm:H^i(G/K^\pm)\to H^i(G/H)$ is onto for $i=2,4$. We will compute the order of $H^4(M)$ from \eqref{ourMVseq} by finding the orders of $\Im(\Delta_3)$ and $\Im(I_4)$.

We will first compute $\Im(I_4)=\ker(\Pi_4:\Z/(h/b_-)\oplus \Z/(h/b_+)\twoheadrightarrow \Z/\delta)$. First notice that $\delta$ divides $h/b_\pm$ because each $\pi_\pm^*: H^4(G/K^\pm)(\simeq\Z/(h/b_\pm)) \to H^4(G/H) (\simeq\Z/\delta)$ is onto. Therefore $\delta|\gcd(h/b_-,h/b_+)$. Conversely, notice that in the explicit example of $\rho$ given in \eqref{rhoeg} we have $\gcd(h/b_-,h/b_+)| \gcd(A,B,D,E)= \delta$. Then since $\delta=|H^4(G/H)|$ is independent of the choice of $\rho$, we have that $\delta=\gcd(h/b_-,h/b_+)$ for any possible $\rho$. It then follows that $\ker(\Pi_4: \Z/(h/b_-)\oplus \Z/(h/b_+) \twoheadrightarrow \Z/\delta)$ is cyclic of order $(\frac{h}{b_-}\cdot \frac{h}{b_+})/ \gcd(\frac{h}{b_-},\frac{h}{b_+}) = \lcm(\frac{h}{b_-},\frac{h}{b_+})$. Therefore $\Im(I_4)\simeq \Z/\lcm(h/b_-,h/b_+)$.

Next consider $\Im(\Delta_3)\simeq H^3(G/H)/\Im(\Pi_3)$. Recall that we have the isomorphisms $H^3(G/K^\pm)\simeq \bE_\infty^{2,1}$ and $H^3(G/H)\simeq E_\infty^{2,1}$. In the ordered bases $\set{u_1u_2}\subset E_2^{0,2}$, $\{v_1u_1, v_1u_2,$ $v_2u_1, v_2u_2\} \subset E_2^{2,1}$ and $\set{v_1v_2}\subset E_2^{4,0}$ we see $\displaystyle{d_2^{2,1} = \left[\begin{array}{cccc} B & E & A & D \end{array}\right] }$, as a matrix, and $\Im(d_2^{0,2})$ is generated by $x:=(-D, A, -E, B)^t$. Further, the image of $\pi_\pm^*$ is generated by $y_\pm:=(q_\pm \hqpm, - q_\pm \hppm, p_\pm \hqpm, - p_\pm \hppm)^t$. Therefore $H^3(G/H)/\Im(\Pi_3) \simeq$ $\ker\displaystyle{\left[\begin{array}{cccc} B & E & A & D \end{array}\right]}/$ $\langle x,y_-,y_+ \rangle$.
Notice further that $y_\pm= \frac{1}{b_\pm}\big( p_\pm q_\pm(D,-A,E,-B)+$ $(q_\pm^2E,-q_\pm^2B,p_\pm^2D,-p_\pm^2A) \big)^t$, after unpacking the definitions of $\hppm,\hqpm$. If we denote $z_\pm:=(q_\pm^2E,-q_\pm^2B,p_\pm^2D,-p_\pm^2A)^t$, then we see that the real span of $x,y_-,y_+$ is the same as that of $x,z_-,z_+$. In particular, if $(p_-^2,q_-^2)=(p_+^2,q_+^2)$ (i.~e.~$(p_-,q_-)=\pom(p_+,-q_+)$) then $x,y_-,y_+$ are linearly dependent over $\R$. Therefore $\Im(\Delta_3)$, and hence $H^4(M)$, are infinite in this case. Otherwise we claim that $x,y_-,y_+$ are linearly independent over $\R$, making $H^4(M)$ finite. To see this notice that when $(p_-^2,q_-^2)\ne(p_+^2,q_+^2)$, we have $\Span_\R\set{z_-,z_+}\ni (E,-B,0,0)^t,(0,0,D,-A)^t$. It is clear that these two vectors and $x$ form a linearly independent set over $\R$, so $\Im(\Delta_3)$ is finite in this case. To find the order of $\Im(\Delta_3)$ we need an integer basis for $\ker\displaystyle{\left[\begin{array}{cccc} B & E & A & D \end{array}\right]}$ and its extension to an integer basis of $\Z^4$. Let $\bA=A/\gcd(A,B)$, $\bB=B/\gcd(A,B)$, $\bD=D/\gcd(D,E)$, $\bE=E/\gcd(D,E)$ and choose integers $\tA,\tB,\tD,\tE$ such that $\bA\tB+\bB\tA=1$ and $\bD\tE+\bE\tD=1$. Then the vectors $w_1:=(\bA,0,-\bB,0)^t$, $w_2:=(\tA,0,\tB,0)^t$, $w_3:=(0,\bD,0,-\bE)^t$ and $w_4:=(0,\tD,0,\tE)^t$ form an integer basis of $\Z^4$, with $w_1,w_3\in \ker \displaystyle{\left[\begin{array}{cccc} B & E & A & D \end{array}\right]}$. We also see that $\mu:=\gcd(A,B)/\delta$ and $\nu:=\gcd(D,E)/\delta$ are relatively prime integers since $\delta=\gcd(A,B,D,E)$, so we can find integers $\tmu,\tnu$ such that $\nu\tmu+\mu\tnu=1$. Then we see that $w_1,\nu w_2-\mu w_4, w_3, \tnu w_2 +\tmu w_4$ is an integer basis for $\Z^4$ with $w_1,\nu w_2-\mu w_4, w_3$ an integer basis for $\ker \displaystyle{\left[\begin{array}{cccc} B & E & A & D \end{array}\right]}$. Therefore we see that the order of $\Im(\Delta_3)\simeq \ker\displaystyle{\left[\begin{array}{cccc} B & E & A & D \end{array}\right]}/ \langle x,y_-,y_+ \rangle$ is equal to $|\det(x,y_-,y_+,\tnu w_2 +\tmu w_4)|$, the determinant of the matrix formed by the vectors $x,y_-,y_+,\tnu w_2 +\tmu w_4$. This gives the explicit value of $|\Im(\Delta_3)|$. It then follow from \eqref{ourMVseq} that the order of $H^4(M)$ is $|\Im(\Delta_3)|\cdot|\Im(I_4)|=|\det(x,y_-,y_+,\tnu w_2 +\tmu w_4)|\cdot \lcm(h/b_-,h/b_+)$.

\begin{conclusion*}
For a manifold $M$ of type $N^7_A$, if $K^-=K^+$ then $M$ has the homology groups of $S^3\times S^2\times S^2$. If $K^-\ne K^+$ we have $H^5(M)=\Z\oplus\Z$ and an exact sequence $0\to \Z/\beta\to H^4(M) \to \Z/\gamma \to 0$ with $\gamma=|ah/b_-b_+|$. If $(p_-,q_-)=\pom(p_+,-q_+)$ then $H^4(M)$ is infinite and hence $\beta=0$. If $(p_-,q_-)\ne \pom(p_+,-q_+)$ then $H^4(M)$ is finite and $\beta=|\det(x,y_-,y_+,\tnu w_2 +\tmu w_4)|\cdot \lcm(h/b_-,h/b_+)/\gamma$.
\end{conclusion*}

\subsection{Actions of type $N^7_B$:}

$$S^3\times S^3 \,\, \supset \,\, \set{(e^{ip\theta},e^{iq\theta})}\cdot H_+, \,\, \set{(e^{j\theta},1)}\cdot H_- \,\, \supset \,\,  H_-\cdot H_+$$
where $\gcd(p,q)=1$, $H_\pm = \Z_{n_\pm} \subset K^\pm_0$, $n_+\le 2$, $4|n_-$ and $p \equiv \pm \frac{n_-}{4} \mod n_-$.

\vspace{.5em}

First notice the following simple but useful fact about this family of actions: $(-1,1)\in H_-$ if and only if $q$ is even. This follows from the facts that $4| n_-$ and $\gcd(p,q)=1$ Furthermore we can assume $q\ne0$ since otherwise the action is a product action.

We will approach this case using the non-primitivity fiber bundle \eqref{nonprim}, with $L=S^3\times S^1$. This bundle has the form $M_L\to M\to S^2$ and we will eventually use the spectral sequence for this bundle to compute the homology of $M$. First we must study $M_L$ and compute its homology groups.

We will start by computing $\pi_1(M_L)$. For this, choose curves $\alpha_\pm:[0,1]\to K^\pm_0$ with $\alpha_\pm(0)=1$ which represent $\pi_1(K^\pm/H)$. By \cite[Prop. 1.7]{thesis} we have $\pi_1(M_L) \simeq \pi_1(L/H) / \grp{\alpha_-,\alpha_+}$. Now let $\delta:[0,1] \to L: t\mapsto (1,e^{2\pi i t})$, a curve which represents $\pi_1(L)$. Further, since $L\to L/H$ is a cover and since $H=H_-\cdot H_+=\grp{\alpha_-(1), \alpha_+(1)}$, it is clear that $\alpha_-$, $\alpha_+$ and $\delta$ generate $\pi_1(L/H)$, when considered as loops in $L/H$. These three elements also commute since $\alpha_-$ and $\delta$ lie in $T^2\subset L$ and $\alpha_+$ can be homotoped into $T^2$. We now divide the computation of $\pi_1(M_L)$ into two cases, depending on whether $q$ is even or odd.

First suppose $q$ is odd and hence $(-1,1)\notin H_-$. If there is some relation $\alpha_-^x\alpha_+^y \delta^z\sim 1$ in $\pi_1(L/H)$ then $\alpha_-^x\alpha_+^y$ must be a loop in $L$. However, since $(-1,1)\notin H_-$, no power of $\alpha_-$ will be a path ending at $(-1,1)$ so $\alpha_-^x$ and $\alpha_+^y$ must both be loops in $L$. It then follows that $\alpha_+^y\subset S^3\times 1$ is contractible in $L$ so we can assume $\alpha_-^x\sim\delta^{-z}$. For $\alpha_-^x$ to be a loop in $L$, $x$ must be a multiple of $n_-=|K^-\cap H|$, say $x=\hat{x}n_-$. It then follows $-z=\hat{x}q$. So the only relations in $\pi_1(L/H)$ are $\alpha_-^{n_-}\sim \delta^q$ and $\alpha_+^{n_+}\sim 0$. Therefore $\pi_1(M_L) \simeq \pi_1(L/H)/\grp{\alpha_-,\alpha_+} \simeq \grp{\delta,\alpha_-, \alpha_+:\delta^q=\alpha_-^{n_-}, \alpha_+^{n_+}=1}_\text{ab}/\grp{\alpha_-, \alpha_+} \simeq \grp{\delta:\delta^q=1}\simeq \Z/q$.

Next suppose $q$ is even and hence $(-1,1)\in H_-$ and $n_+=2$. This time $\alpha_-^{n_-/2}$ and $\alpha_+$ are both paths from $(1,1)$ to $(-1,1)$. Hence we have the additional relation $\alpha_-^{n_-/2}\sim \alpha_+\delta^{q/2}$. As in the previous case $\pi_1(M_L) \simeq \grp{\delta,\alpha_-, \alpha_+:\alpha_-^{n_-/2}\sim \alpha_+\delta^{q/2}, \alpha_+^2=1}_\text{ab}/\grp{\alpha_-, \alpha_+} \simeq \grp{\delta:\delta^{q/2}=1}\simeq \Z/(q/2)$. Notice that whether $q$ is even or odd we have the following unifying equation: $\pi_1(M_L)=\Z/(q/\gcd(q,2))$.

Now, in order to compute the homology groups of $M_L$, we need the homology groups of $L/K^\pm$. First $L/K^-= S^3\times S^1/\set{(e^{ip\theta},e^{iq\theta})}\cdot \set{((-1)^{n_++1},1)}$. We see $S^3$ acts on $L/K^-$ transitively since $q\ne 0$. Since $\gcd(p,q)=1$, the isotropy is $S^3\times 1\cap \set{(e^{ip\theta},e^{iq\theta})}\cdot \set{((-1)^{n_++1},1)}\simeq \Z/\lcm(n_+,q)$, so $L/K^-\simeq S^3/\Z_{\lcm(n_+,q)}$, a lens space.

For $L/K^+$ we can write $L/K^+=S^3\times S^1/\set{(e^{j\theta},1)}\cdot \set{(\xi_{n_-}^p,\xi_{n_-}^q)}$ where $\xi_{n_-}$ is a primitive $n_-^\text{th}$ root of unity. Since $p_- \equiv \pm \frac{n_-}{4} \mod n_-$, we know $\xi_{n_-}^p=\pm i$. Then we notice $S^3/\set{e^{j\theta}}\simeq S^2$ and left multiplication by $\pm i$ on $S^3/\set{e^{j\theta}}\simeq S^2$ acts as $-\Id$. Therefore $L/K^+\simeq S^2\times S^1/\grp{(-\Id,\xi_{n_-}^q)}$. If $b:=n_-/\gcd(q,n_-)$, the order of $\xi_{n_-}^q$, is odd then $\grp{(-\Id,\xi_{n_-}^q)}=\grp{(-\Id,1), (\Id,\xi_{n_-}^q)}$ so $L/K^+\simeq \RP^2\times S^1$ in this case. If $b$ is even, then $L/K^+\simeq S^2\times S^1/\grp{(-\Id,-1)}$ and this space can be modeled as $S^2\times [0,1]/\sim$ where $(x,0)\sim (-x,1)$. Using the Mayer-Vietoris sequence for the decomposition $S^2\times [0,1/2]$ and $S^2\times [1/2,1]$ we see that the nontrivial cohomology groups for $L/K^+$ are $H^0(L/K^+)\simeq H^1(L/K^+)\simeq \Z$ and $H^3(L/K^+)\simeq \Z_2$, in the case that $b$ is even. When $b$ is odd, we see $L/K^+\simeq \RP^2\times S^1$ has these same cohomology groups with one additional nontrivial group $H^2(L/K^+)\simeq \Z_2$.

We will compute the cohomology groups of $M_L$ using the modified long exact sequence of the pair \eqref{lespair}, but we must first show that $K^-/H\to L/H\to L/K^-$ is orientable as a sphere bundle, in the sense of \cite[pg. 442]{Hatcher}. For this, consider the diagram of bundles \eqref{showbundleorient}, with $F=L$ and $T=T^2$.
\begin{equation}\label{showbundleorient}
\displaystyle{\begin{array}{ccccc}
      K^-/H & \to & F/H & \to & F/K^-\\
      \uparrow & & \uparrow & & \uparrow \\
      K^-/H & \to & T/H & \to & T/K^-
 \end{array}  }
\end{equation}
The bottom bundle is equivalent to the standard product bundle $S^1\to T^2\to S^1$ which is orientable. Notice also that the last vertical map induces a map $\pi_1(T/K^-)\to\pi_1(L/K^-)$, which is onto. Hence any loop in $L/K^-$ can be homotoped to a loop $\gamma$ in $T/K^-$. Then the homotopy of $K^-/H$ above $\gamma$ can be chosen as the product homotopy in $T/H$. It then follow that this induces the identity map on the fiber $K^-/H\to K^-/H$. Hence the top bundle in \eqref{showbundleorient} is also orientable as a sphere bundle.

Hence we have the sequence \eqref{lespair} with $L$ in place of $G$ and $M_L$ in place of $M$. This diagram clearly gives $H^5(M_L)\simeq \Z$ so that $M_L$ is orientable. Then, using the fact that $H^4(M_L)\simeq H_1(M_L)\simeq \pi_1(M_L)\simeq \Z/(q/\gcd(q,2))$ the sequence \eqref{lespair} shows $H^3(M_L)\simeq\Z_2$ if $n_+=1$ and $q$ is odd and $H^3(M_L)=0$ if $n_+=2$ or $q$ is even. It then follows $H^2(M_L) \simeq \Z/(q/\gcd(q,2))$ and $H^1(M_L)=0$.

We are now in a position to consider the spectral sequence for the non-primitivity bundle $M_L\to M\to S^2$. The first page of this spectral sequence takes the form $E^{p,q}_2\simeq H^p(S^2)\otimes H^q(M_L)$. It is clear that only the maps on the first page $d_2^{p,q}:E^{p,q}_2\to E^{p+2, q-1}_2$ are potentially nontrivial, so $E_3=E_\infty$. Notice first that $E_3^{0,5}=\ker\left(d_2^{0,5}:\Z\to \Z/(q/\gcd(q,2))\right)\simeq \Z$. We now consider two cases depending on $H^3(M_L)$. If $n_+=2$ or $q$ is even then $H^3(M_L)=0$, and so $E_2^{0,3}=E_2^{2,3}=0$. It then follows that $H^5(M)\simeq E_3^{0,5}\simeq \Z$ since this is the only nontrivial group along that diagonal. Further $d_2^{0,4}=d_2^{0,3}=0$ so $E_2^{0,4}=E_3^{0,4}$ and $E_2^{2,2}=E_3^{2,2}$. Therefore we have the filtration $H^4(M)=D^{0,4}\supset D^{1,3}= D^{2,2}\supset D^{3,1}= D^{4,0} =0$ with $D^{0,4}/D^{1,3}\simeq E_3^{0,4}\simeq \Z/(q/\gcd(q,2))$ and $D^{2,2}/D^{3,1}\simeq E_3^{2,2}\simeq \Z/(q/\gcd(q,2))$. Hence $H^4(M)$ fits into the short exact sequence $0\to\Z/(q/\gcd(q,2)) \to H^4(M) \to \Z/(q/\gcd(q,2)) \to 0$. Next, if $n_+=1$ and $q$ is odd then $H^3(M_L)\simeq \Z_2$. We then see that $d_2^{0,4}:\Z/q\to \Z/2$ and $d_2^{0,3}:\Z/2\to \Z/q$ are trivial since $q$ is odd. So $E_3^{0,4}=E_2^{0,4}\simeq \Z/q$, $E_3^{2,2}=E_2^{2,2}\simeq \Z/q$ and $E_3^{2,3}=E_2^{2,3}\simeq \Z/2$. Then, as before, we have the short exact sequence $0\to\Z/q \to H^4(M) \to \Z/q \to 0$. Similarly we see $H^5(M)\simeq \Z\oplus \Z_2$.

\begin{conclusion*}
A manifold $M$ of type $N^7_B$ has $H^5(M)\simeq \Z\oplus (\Z/\alpha)$ where $\alpha=1$ if $n_+=2$ or $q$ is even; and $\alpha =2$ if $n_+=1$ and $q$ is odd. Further there is an exact sequence $0\to\Z/(q/\gcd(q,2)) \to H^4(M) \to \Z/(q/\gcd(q,2)) \to 0$.
\end{conclusion*}

\subsection{Actions of type $N^7_C$:}

$$S^3\times S^3 \,\, \supset \,\, \set{(e^{ip\theta},e^{iq\theta})}, \,\, S^3\times \Z_n \,\, \supset \,\,  \Z_n$$
where $\gcd(p,q)=1$, $\gcd(q,n)=1$ and we can assume $q\ne 0$ otherwise this is a product action.

\vspace{.5em}

This case will require several steps. Notice first that $G/K^-\simeq S^3\times S^2$ \cite[Prop. 2.3]{WZ} and $G/K^+\simeq S^3/\Z_n$ is a lens space. Since $G/K^-$ is simply connected we get the long exact sequence of the pair \eqref{lespair}. This sequence easily shows that $H^5(M)\simeq \Z$ and $H^4(M)$ is cyclic. It only remains to find the order of $H^4(M)$.

For this, we will use the non-primitivity fiber bundle \eqref{nonprim} with $L=S^3\times S^1$. This takes the form $M_L\to M\to S^2$ where $M_L$ is given by the diagram $S^3\times S^1\supset \set{(e^{ip\theta},e^{iq\theta})}, S^3\times \Z_n\supset\Z_n$. Clearly $L/K^+\simeq S^1$ and we claim $L/K^-\simeq S^3/\Z_q$. To see this note $S^3$ acts on $L/K^-$ in the obvious way. This action is transitive since $q\ne 0$ and the isotropy group is $K^-\cap(S^3\times 1)=\Z_q \times 1$, since $\gcd(p,q)=1$, and this proves the claim.

We now claim that $K^-/H\to L/H\to L/K^-$ is orientable as a sphere bundle. To see this, notice we have the diagram of sphere bundles \eqref{showbundleorient}, with $F=L$ and $T=S^1\times S^1\subset L$. By the same argument given below \eqref{showbundleorient}, we see $L/H\to L/K^-$ is orientable.

Therefore we get the long exact sequence for the pair \eqref{lespair} with $L$ in place of $G$, and $M_L$ in place of $M$. This sequence clearly gives $H^5(M_L)\simeq \Z$, $H^4(M_L)\simeq \Z_q$ and $H^3(M_L)=0$. It then follows that $M_L^5$ is orientable and hence $H^2(M_L)\simeq \Z_q$ and $H^1(M_L)=0$.

We can then plug this information into the spectral sequence for the bundle $M_L\to M\to S^2$. We see that $E_\infty^{0,4}\simeq E_2^{0,4}\simeq H^4(M_L)\simeq \Z_q$, $E_\infty^{2,2}\simeq E_2^{2,2}\simeq H^2(S^2)\otimes H^2(M_L)\simeq \Z_q$ and $E_\infty^{i,4-i}\simeq E_2^{i,4-i}=0$ for $i\ne0,2$. It follows that $H^4(M)$ as a subgroup $\Z_q\subset H^4(M)$ with $H^4(M)/\Z_q\simeq \Z_q$. In particular $H^4(M)$ has order $q^2$. Since we already showed $H^4(M)$ is cyclic, it follows $H^4(M)\simeq \Z/q^2$.

\begin{conclusion*}
If $M$ is a manifold of type $N^7_C$, then $H^5(M)\simeq \Z$ and $H^4(M)\simeq \Z/q^2$.
\end{conclusion*}

\subsection{Actions of type $N^7_D$:}

$$S^3\!\!\times\! S^3\!\!\times\! S^1 \supset  \set{(z^mw^{\mu p}\!,z^nw^{\nu p}\!,w)},   \set{(z^mw^{\mu p}\!,z^nw^{\nu p}\!,w)}   \supset  H_0\cdot\Z_a $$
where $H_0=\set{(z^m,z^n,1)}$, $m\nu-n\mu=1$ and $\Z_a\subset \set{(w^{\mu p}\!,w^{\nu p}\!,w)}$.

\vspace{.5em}

We first claim that $G/K^\pm\simeq S^3\times S^2$. To see this, note that $S^3\times S^3$ acts on $G/K^\pm$ by multiplication on the first two components in the natural way. It is clear this action is transitive and that the isotropy is $K^\pm\cap (S^3\times S^3\times 1)=\set{(z^m,z^n,1)}$. Hence $G/K^\pm\simeq S^3\times S^3/\set{(z^m,z^n)}$ which is know to be $S^3\times S^2$, as before. In particular $G/K^\pm$ is simply connected so we get the long exact sequence of the pair \eqref{lespair}. One segment of this sequence is $0\to\Z\to H^5(M)\to\Z\to 0$ and hence $H^5(M)\simeq \Z^2$.

We also have the non-primitivity fiber bundle for $L=K^\pm$ \eqref{nonprim}, which takes the form $S^2\to M\to S^3\times S^2$ in this case. The Gysin sequence of this bundle contains the segment $0\to H^4(M)\to\Z\to\Z\to \Z^2\to \Z \to 0$. Hence $H^4(M)\simeq \Z$. In fact it follows that the Euler class of this bundle is trivial.

\begin{conclusion*}
A manifold $M$ of type $N^7_D$ has the same homology groups as $S^3\times S^2\times S^2$. 
\end{conclusion*}

\subsection{Actions of type $N^7_E$:}

$$S^3\!\!\times\! S^3\!\!\times\! S^1 \supset  \set{\!(z^mw^{\mu p_-}\!,z^nw^{\nu p_-}\!,w^{q_-})\!}\! H,   \set{\!(z^mw^{\mu p_+}\!,z^nw^{\nu p_+}\!,w^{q_+})\!}\! H  \supset H$$
where $\gcd(p_\pm,q_\pm)=1$, $H=H_-\cdot H_+$, $H_0=\set{(z^m,z^n,1)}$, $K^-\ne K^+$, $m\nu-n\mu=1$, $\gcd(q_-,q_+,d)=1$ where $d$ is the index of $H\cap K^-_0 \cap K^+_0$ in $K^-_0 \cap K^+_0$.

\vspace{.5em}

This case will be very similar to the case of family $N^7_A$. As a result we will skip many of the details here and, instead, indicate the significant modifications from the family $N^7_A$. Note first that $L=T^3\subset G$ contains $K^\pm$ and $H$. Now let $\tL:=\set{(w_1^{\mu},w_1^{\nu},w_2)}$ and notice that $L$ is the direct product $L= H_0\cdot \tL$, since $m\nu-n\mu=1$. Also denote $\tK^\pm:=K^\pm\cap\tL$, $\tH:=H\cap \tL$ and $\tH_\pm:=H_\pm\cap \tL$. Now define $h:=|\tH|=|H/H_0|$ and $b_\pm:=|\tH_\pm|=|(K^\pm_0\cap H)/H_0|$ so that $K^\pm$ and $\tK^\pm$ have $h/b_\pm$ connected components. Finally define $a:=q_+p_--q_-p_+$ and notice that $\tK^-_0\cap\tK^+_0$ has order $|a|$.

Now consider the non-primitivity fiber bundle \eqref{nonprim}, $M_L\to M\to G/L$, where $G/L\simeq S^2\times S^2$. The group diagram for $M_L$ reduces to $\tL\supset \tK^-,\tK^+\supset \tH$, so $M_L\simeq S^3/\Z_r$ where $r=|ha/b_-b_+|$, just as in case $N^7_A$. Then, by the argument given for $N^7_A$, we see that $H^5(M)\simeq \Z\oplus\Z$ and $0\to \Z/\beta \to H^4(M)\to \Z/\gamma\to 0$ is short exact where $\gamma=r=|ha/b_-b_+|$. Again, it only remains to compute the value of $\beta$.

Again we have the diagram \eqref{pidiagram}, this time with $L=T^3$, etc. As before, denote $\pi_\pm: G/H\to G/K^\pm$ as the natural projection, and let $(E_i^{k,l},d_i^{k,l})$ and $(\bE_i^{k,l},\bd_i^{k,l})$ be the spectral sequences for the top and the bottom bundles of \eqref{pidiagram} respectively. To choose our generators of $H^1(L/H)$, label the map $\sigma:L\to \tL:(z^m w_1^{\mu},z^n w_1^{\nu},w_2)\mapsto (w_1^{\mu},w_1^{\nu},w_2)$, so $\ker(\sigma)=H_0$, and choose a homomorphism $\rho:\tL\to T^2:(w_1^{\mu},w_1^{\nu},w_2)\mapsto (w_1^Aw_2^B,w_1^Dw_2^E)$ with $\ker(\rho)=\tH$. For example, we may chose $A,B,D,E$ precisely as in \eqref{rhoeg}, just as before. Then $\rho\circ\sigma: L\to T^2$ has $\ker(\rho\circ\sigma)=H$, so we have an isomorphism $f:L/H\to T^2:\ell H\to \rho(\sigma(\ell))$. Define $u_i:=f^*(\pr_i^*(\bbone))$, $i=1,2$, as before, for some generator $\bbone\in H^1(S^1)$. Clearly $u_1,u_2$ generate $H^1(L/H)$. Finally, just as with $N^7_A$, there are generators $u_\pm\in H^1(L/K^\pm)\simeq \Z$ with $\pi_\pm^*(u_\pm)=\hqpm u_1-\hppm u_2$, where $\hppm=(Ap_\pm+Bq_\pm)/b_\pm$ and $\hqpm=(Dp_\pm+Eq_\pm)/b_\pm$ are relatively prime integers.

We will use the spectral sequence $(\bE_i^{k,l},\bd_i^{k,l})$ to compute $H^i(G/K^\pm)$, but first we must determine the map $\bd_2^{0,1}$. For this consider the bundle diagram \eqref{hopfdiagram} with $F=K^\pm$, $L=T^3$, etc., and denote the spectral sequence of the top bundle as $(\tE_i^{k,l},\td_i^{k,l})$. As before, fix generators $\bbone\in H^1(S^1)$ and $\bbone'\in H^2(S^2)$ and define the standard product generators $\tu_i:=\pr_i^*(\bbone)\in H^1(L)$, $i=1,2,3$, and $v_i:=\pr_i^*(\bbone')\in H^2(G/L\simeq S^2\times S^2)$, $i=1,2$. With the right choice of sign on $\bbone'$ we have $\td_2^{0,1}(\tu_1)=v_1$, $\td_2^{0,1}(\tu_2)=v_2$ and $\td_2^{0,1}(\tu_3)=0$.
As in the case $N^7_A$, we see the induced map $\wp_\pm^*:H^1(L/K^\pm)\to H^1(L)$ takes a generator to $\pom\frac{h}{b_\pm}(q_\pm n\tu_1-q_\pm m\tu_2+p_\pm \tu_3)$.
Hence $\bd_2^{0,1}(u_\pm)=\td_2^{0,1}(\wp_\pm^*(u_\pm))=\pom\frac{hq_\pm}{b_\pm}(n v_1- m v_2)$, using the commutativity of the spectral sequence diagram.

We will now compute $H^i(G/K^\pm)$ in two cases, depending on whether or not $q_\pm=0$. If $q_\pm=0$ then $\bd_2^{0,1}=0$ and it is clear from $\bE_2^{k,l}$ that $H^1(G/K^\pm)\simeq \Z$, $H^2(G/K^\pm)\simeq \Z\oplus\Z$, $H^3(G/K^\pm)\simeq \Z\oplus\Z$, $H^4(G/K^\pm)\simeq \Z$ and $H^5(G/K^\pm)\simeq \Z$.
If $q_\pm\ne 0$ then $\bd_2^{0,1}\ne 0$ and we can compute $\bd_2^{2,1}(v_1u_\pm)= \mop\frac{hmq_\pm}{b_\pm}v_1v_2$ and $\bd_2^{2,1}(v_2u_\pm)= \pom\frac{hnq_\pm}{b_\pm}v_1v_2$, where $v_1u_\pm,v_2u_\pm\in \bE_2^{2,1}$ generate. Hence $\ker(\bd_2^{2,1})= \grp{nv_1u_\pm+mv_2u_\pm}\simeq \Z$, since $m$ and $n$ are coprime. We then see from $\bE_2^{k,l}$ that $H^1(G/K^\pm)=0$, $H^2(G/K^\pm)\simeq \grp{v_1,v_2}_\ab / \grp{\frac{hq_\pm}{b_\pm}(n v_1- m v_2)} \simeq \Z \oplus \big(\Z/(hq_\pm/b_\pm)\big)$, $H^3(G/K^\pm)\simeq \grp{nv_1u_\pm+mv_2u_\pm}\simeq \Z$, $H^4(G/K^\pm)\simeq \grp{v_1v_2}/ \grp{\frac{hnq_\pm}{b_\pm}v_1v_2, \frac{hmq_\pm}{b_\pm}v_1v_2}= \grp{v_1v_2}/ \grp{\frac{hq_\pm}{b_\pm}v_1v_2}\simeq \Z/(hq_\pm/b_\pm)$ and $H^5(G/K^\pm)\simeq \Z$.

To repeat this process for $G/H$, consider the diagram \eqref{hopfdiagram} with $F=H$, $L=T^3$, etc., and denote $\wp:L\to L/H$. We see first that $\wp^*(u_i)=\wp^*(f^*(\pr_i^*(\bbone)))= (\rho\circ\sigma)^*(\pr_i^*(\bbone))$, since $\rho(\sigma(\ell))= f(\wp(\ell))$. We can then compute $\wp^*(u_1)= -An\tu_1+Am\tu_2+B\tu_3$ and $\wp^*(u_2)= -Dn\tu_1+Dm\tu_2+E\tu_3$. Hence $d_2^{0,1}(u_1) = \td_2^{0,1}(\wp^*(u_1)) = A(-nv_1+mv_2)$ and $d_2^{0,1}(u_2) = D(-nv_1+mv_2)$. We then have $d_2^{0,2}(u_1u_2)= Dnv_1u_1 -Anv_1u_2 -Dmv_2u_1 +Amv_2u_2$. We also get $d_2^{2,1}(v_1u_1)=Amv_1v_2$, $d_2^{2,1}(v_1u_2)=Dmv_1v_2$, $d_2^{2,1}(v_2u_1)=-Anv_1v_2$ and $d_2^{2,1}(v_2u_2)=-Dnv_1v_2$. Finally we compute $d_2^{2,2}(v_1u_1u_2)=-Dmv_1v_2u_1+Amv_1v_2u_2$ and $d_2^{2,2}(v_2u_1u_2)=Dnv_1v_2u_1-Anv_1v_2u_2$.

We can now compute $H^i(G/H)$ using $E_2^{k,l}$. We easily get $H^1(G/H)\simeq \Z$, $H^2(G/H)\simeq \grp{v_1,v_2}_\ab/\grp{A(-nv_1+mv_2), D(-nv_1+mv_2)}\simeq \Z\oplus(\Z/\ell)$ where $\ell=\gcd(A,D)$, $H^3(G/H)\simeq \ker(d_2^{1,2})/\Im(d_2^{0,2})\simeq \Z\oplus\Z\oplus(\Z/\ell)$ by Poincar\'{e} duality, $H^4(G/H)\simeq \Z\oplus(\Z/\ell)$, $H^5(G/H)\simeq \Z\oplus(\Z/\ell)$ and $H^6(G/H)\simeq \Z$.

Now consider the Mayer-Vietoris sequence \eqref{MVseq} for $M$. It will take the form \eqref{ourMVseq} and we will again compute the order of $H^4(M)$ by computing the orders of $\Im(I_4)$ and $\Im(\Delta_3)$. As in case $N^7_A$ it is clear that $\Delta_2=0$ since $\pi_\pm: H^2(G/K^\pm)\to H^2(G/H)$ are both onto, for the same reason as before. There will be two cases to consider, depending on whether or not $q_-q_+=0$. First, suppose $q_-q_+=0$. Since $K^-\ne K^+$ by assumption it follows that only one of $q_-$ or $q_+$ can be zero. It would then follow that $H^3(G/K^-)\oplus H^3(G/K^+)\simeq \Z^3$ and $H^3(G/H)\simeq \Z\oplus \Z \oplus(\Z/\ell)$, with $\ell\ne 0$. Hence, from \eqref{ourMVseq}, $H^3(M)$ would be infinite. Therefore $H^4(M^7)$ would be infinite as well, and so $\beta=0$ in this case.

Now assume that $q_-q_+\ne 0$. First we will find the order of $\Im(I_4)$. Notice from the commutativity of the spectral sequences from \eqref{pidiagram}, that $\pi_\pm: H^4(G/K^\pm)\to H^4(G/H)$ maps $H^4(G/K^\pm)\simeq \grp{v_1v_2}/ \grp{\frac{hq_\pm}{b_\pm}v_1v_2}\simeq \Z/(hq_\pm/b_\pm)$ onto the torsion subgroup of $H^4(G/H)$, $\grp{v_1v_2}/ \grp{\ell v_1v_2}\simeq \Z/\ell$. It follow that the order of $\Im(I_4)=\ker(\Pi_4)$ is $(\frac{hq_-}{b_-}\frac{hq_+}{b_+})/\ell=h^2q_-q_+/b_-b_+\ell$. In fact, one can show that $\ell=\gcd(hq_-/b_-,hq_+/b_+)$ so that $\Im(I_4)$ is cyclic, but this will not be important for our purposes.

Finally we must find the order of $\Im(\Delta_3)\simeq H^3(G/H)/\Im(\Pi_3)$, from \eqref{ourMVseq}. We will do this in precisely the same manner as we did with $N^7_A$. In the natural bases, we see $d_2^{2,1}=\displaystyle{\left[\begin{array}{cccc} Am & Dm & -An & -Dn \end{array}\right]}$ and $\Im(d_2^{0,2})$ is generated by $x:=(Dn, -An, -Dm, Am)^t$. Furthermore, from the commutativity of the spectral sequence diagram we get that $\Im(\pi_\pm^*:\bE_3^{2,1}\to E_3^{2,1})=\pi_\pm^*(\ker(\bd_2^{2,1}))$. We recall that $\ker(\bd_2^{2,1})= \grp{nv_1u_\pm+mv_2u_\pm}$ and $\pi_\pm^*(u_\pm)=\hqpm u_1-\hppm u_2$. Hence $\Im(\pi_\pm^*)\subset E_3^{2,1}$ is generated by $y_\pm:= (n\hqpm,-n\hppm,m\hqpm,-m\hppm)^t$. Therefore $\Im(\Delta_3)\simeq H^3(G/H)/\Im(\Pi_3) \simeq \ker\displaystyle{\left[\begin{array}{cccc} Am & Dm & -An & -Dn \end{array}\right]} / \grp{x,y_-,y_+}$. This will be finite if and only if $\set{x,y_-,y_+}$ is linearly independent over $\R$. Notice that by the definition of $(\hppm,\hqpm)$, the pairs $\set{(\hpm,\hqm),(\hpp,\hqp)}$ are linearly independent. Also $y_\pm= \hqpm(n,0,m,0)^t- \hppm(0,n,0,m)^t$ and $x:=D(n, 0, -m, 0)^t+A(0, -n, 0, m)^t$. We see then that $\set{x,y_-,y_+}$ is linearly independent over $\R$ if and only if $mn\ne 0$. In particular if $mn=0$ then $\Im(\Delta_3)$ is infinite. Therefore $H^4(M)$ is infinite and hence $\beta=0$ in this case. Otherwise $H^4(M)$ is finite. To find its order, notice that $w_1:=(n,0,m,0)^t$, $w_2:=(0,n,0,m)^t$, $w_3:=(\nu,0,\mu,0)^t$ and $w_4:=(0,\nu,0,\mu)$ form an integer basis for $\Z^4$ with $w_1,w_2\in \ker \displaystyle{ \left[\begin{array}{cccc} Am & Dm & -An & -Dn \end{array}\right] }$. Let $\bA:=A/\ell$ and $\bD:=D/\ell$ then $\gcd(\bA,\bD)=1$ since $\ell=\gcd(A,D)$, so we can choose integers $\zeta,\eta$ with $\bA\zeta+\bD\eta=1$. Then $\set{w_1,w_2,\bD w_3-\bA w_4,\zeta w_3+\eta w_4}$ is also an integer basis for $\Z^4$, with $\set{w_1,w_2,\bD w_3-\bA w_4}$ an integer basis for $\ker \displaystyle{ \left[\begin{array}{cccc} Am & Dm & -An & -Dn \end{array}\right] }$. Therefore the order of $\Im(\Delta_3)\simeq \ker\displaystyle{ \left[\begin{array}{cccc} Am & Dm & -An & -Dn \end{array}\right]} / \grp{x,y_-,y_+}$ is given by $|\det(x,y_-,y_+,\zeta w_3+\eta w_4)|$, the absolute value of the determinant of the matrix built from the four vectors $x,y_-,y_+,\zeta w_3+\eta w_4\in \Z^4$. Therefore $|H^4(M)|=|\Im(\Delta_3)|\cdot|\Im(I_4)|= |\det(x,y_-,y_+,\zeta w_3+\eta w_4)| h^2q_-q_+/b_-b_+\ell$.

\begin{conclusion*}
If $M$ is a manifold of type $N^7_E$, then $H^5(M)\simeq \Z\oplus \Z$ and there is an exact sequence $0\to \Z/\beta \to H^4(M)\to \Z/\gamma\to 0$ with $\gamma=|ha/b_-b_+|$. Furthermore, $\beta=0$ when $q_-q_+mn=0$, otherwise $\beta$ is finite and $\beta=|H^4(M)|/\gamma=|\det(x,y_-,y_+,\zeta w_3+\eta w_4)| hq_-q_+/\ell |a|$.
\end{conclusion*}

\subsection{Actions of type $N^7_F$:}

$$S^3 \! \times \! S^3 \! \times \! S^1 \,\, \supset \,\, \set{(e^{ip\phi}e^{ia\theta},e^{i\phi},e^{i\theta})}, \,\, S^3\! \times \! S^1\! \times\! \Z_n  \,\, \supset \,\, \set{(e^{ip\phi},e^{i\phi},1)}\cdot \Z_n$$
where $\Z_n\subset \set{(e^{ia\theta},1,e^{i\theta})}$.

\vspace{.5em}

Consider the non-primitivity diagram \eqref{nonprim} with $L=S^3\times S^1\times S^1$. The group diagram for $M_L$ can be reduced to $S^3\times S^1\supset \set{(w^a,w)}, S^3\times \Z_n \supset \Z_n$ using \cite[Prop. 1.12]{thesis}. We then recognize this as an action on $S^5$. Hence \eqref{nonprim} becomes $S^5\to M\to S^2$. The Gysin sequence for this bundle clearly gives the cohomology groups of $M$.

\begin{conclusion*}
A manifold $M$ of type $N^7_F$ has the same homology groups as $S^5\times S^2$.
\end{conclusion*}


\subsection{The action $N^7_G$:}

$$\SU(3) \,\, \supset \,\, \S(\U(1)\U(2)), \,\, \S(\U(1)\U(2)) \,\, \supset \,\, T^2.$$

\vspace{.5em}

We will handle this case in two steps. First notice that the non-primitivity fiber bundle \eqref{nonprim} for $L=\S(\U(1)\U(2))$ is $S^3\to M\to \CP^2$. The Gysin sequence for the bundle easily shows $H^5(M)\simeq \Z$. Now, since $G/K^\pm\simeq \CP^2$ is simply connected the bundle $G/H\to G/K^\pm$ is an orientable sphere bundle. So we have the long exact sequence for the pair $(M,B_+)$, \eqref{lespair}. One segment of this sequence is $0\to H^4(M)\to \Z\to\Z\to \Z\to 0$ which shows $\H^4(M)\simeq \Z$. In fact it follows that the Euler class of the bundle $S^3\to M\to \CP^2$ is trivial.

\begin{conclusion*}
The manifold $M$ of type $N^7_G$ has the homology groups of $\CP^2\times S^3$. 
\end{conclusion*}

\subsection{Actions of type $N^7_H$:}

$$\SU(3)\times S^1 \,\, \supset \,\, \set{(\beta(m_-\theta),e^{in_-\theta})}\cdot H, \,\, \set{(\beta(m_+\theta),e^{in_+\theta})}\cdot H  \,\, \supset \,\, H$$
where $\gcd(m_\pm,n_\pm)=1$, $H_0 = \SU(1)\SU(2)\times 1$, $H=H_-\cdot H_+$, $K^-\ne K^+$, $\beta(\theta)= \diag(e^{-i\theta},e^{i\theta},1)$, and $\gcd(n_-,n_+,d)=1$ where $d$ is the index of $H\cap K^-_0 \cap K^+_0$ in $K^-_0 \cap K^+_0$.

\vspace{.5em}

We will handle this case in several steps. Start by taking $L=\S(\U(1)\U(2))\times S^1$. This gives the non-primitivity bundle \eqref{nonprim} $M_L\to M\to \CP^2$. Since $\SU(1)\SU(2)\times 1$ is normal in all of $L$, $K^\pm$ and $H$, we see the diagram for $M_L$ has the effective form $T^2\supset S^1_-, S^1_+\supset 1$, where $S^1_-$ and $S^1_+$ are distinct circle subgroups of $T^2$. It then follows that $M_L\simeq S^3/\Z_r$ is a lens space for some $r\in\Z_+$ \cite[Sec. 7.2]{thesis}.

Now consider the long exact sequence of homotopy groups for the bundle $S^3/\Z_r \to M\to \CP^2$. This sequence clearly shows $\pi_2(M)\simeq \Z$ and $\pi_3(M)\simeq \Z$. Hurewicz Theorem then implies that $h:\pi_2(M)\to H_2(M)$ is an isomorphism and $h:\pi_3(M)\to H_3(M)$ is onto. In particular $H_2(M)\simeq H^5(M)\simeq \Z$ and $H_3(M)\simeq H^4(M)$ is cyclic, so it only remains to find the order of $H^4(M)$.

We will eventually find the order of $H^4(M)$ by using \eqref{lespair} but first we need to compute the homology groups of $G/K^\pm$. To do this consider the fiber bundle $L/K^\pm\to G/K^\pm\to G/L$, where $L/K^\pm\simeq S^1$ and $G/L\simeq \CP^2$. The Gysin sequence for this bundle shows that $H^5(G/K^\pm)\simeq \Z$, which means $G/K^\pm$ is orientable. This sequence also contains the segment $0\to H^3(G/K^\pm)\to \Z\to \Z\to H^4(G/K^\pm)\to 0$. If the middle map here is multiplication by $k_\pm\in \Z$ then $\Z/k_\pm\simeq H^4(G/K^\pm)\simeq H_1(G/K^\pm)$ and $H^3(G/K^\pm)$ is trivial if $k_\pm\ne0$ or $\Z$ if $k_\pm=0$.

We will find the value of $k_\pm$ by computing $\pi_1(G/K^\pm)$. For this, we first claim $\pi_1(G/K^\pm_0)$ is infinite if $n_\pm=0$ and has order $n_\pm$ if $n_\pm\ne0$. The case $n_\pm=0$ is clear, so suppose $n_\pm\ne0$. Notice that $\SU(3)$ acts on $G/K^\pm_0$ in the natural way. Since $n_\pm\ne0$ this action is transitive and since $\gcd(m_\pm,n_\pm)=1$ the isotropy subgroup is $\SU(1)\SU(2)\cdot \Z_{n_\pm}$. Hence $G/K^\pm_0\simeq \SU(3)/(\SU(1)\SU(2)\cdot \Z_{n_\pm})$, so $\pi_1(G/K^\pm_0)\simeq \Z_{n_\pm}$ in this case, as we claimed. Now consider the covering space bundle $K^\pm/K^\pm_0\to G/K^\pm_0\to G/K^\pm$. If $n_\pm=0$ we see that $\pi_1(G/K^\pm)$ is infinite. If $n_\pm\ne0$ this shows that the order of $\pi_1(G/K^\pm)$ is $n_\pm a_\pm$ where $a_\pm= \left|K^\pm/K^\pm_0\right|$ is the number of connected components of $K^\pm$. From the homotopy exact sequence for $L/K^\pm\to G/K^\pm\to G/L$ we see that $\pi_1(G/K^\pm)$ is cyclic. Putting all this together, we conclude $\pi_1(G/K^\pm)\simeq \Z/(n_\pm a_\pm)$, where we remember $n_\pm$ might be zero here. In particular $k_\pm=n_\pm a_\pm$.

Before we can use \eqref{lespair} we need to show that the sphere bundle $K^\pm/H\to G/H\to G/K^\pm$ is orientable. As before, notice we have the diagram of sphere bundles \eqref{showbundleorient} with $F=G$ and $T=L$. Again the bottom bundle is equivalent to the product bundle $S^1\to T^2\to S^1$, and $\pi_1(L/K^\pm)\to \pi_1(G/K^\pm)$ is onto. The argument given in that case also works in this case. Hence this bundle is orientable as a sphere bundle and so we get the long exact sequence of the pair \eqref{lespair}. If $n_-$ and $n_+$ are both zero then $K^-$ would equal $K^+$ which is impossible by assumption. So we can assume $n_-\ne0$ and hence $k_-\ne0$. Then in \eqref{lespair} we have the segment  $H^3(G/K^+)\to \Z/k_- \to H^4(M) \to \Z/k_+\to 0$. If $k_+=0$ then $H^4(M)$ is infinite. If $k_+\ne0$ then $H^3(G/K^+)=0$ and we see $H^4(M)$ has order $k_-k_+$. Since we already showed $H^4(M)$ is cyclic, we can conclude $H^4(M)\simeq \Z/(k_-k_+)=\Z/(a_-a_+n_-n_+)$, where we remember $n_-n_+$ might be zero.

\begin{conclusion*}
If $M$ is a manifold of type $N^7_H$, then $H^5(M)\simeq \Z$ and $H^4(M)\simeq \Z/(a_-a_+n_-n_+)$ where $a_\pm= \left|K^\pm/K^\pm_0\right|$. In particular if $n_-n_+=0$, $M$ has the homology of $\CP^2\times S^3$.
\end{conclusion*}

\subsection{The action $N^7_I$:}\label{N^7_I}

$$\Sp(2)\,\, \supset \,\, \Sp(1)\Sp(1), \,\, \Sp(1)\Sp(1) \,\, \supset \,\, \Sp(1)\SO(2).$$

\vspace{.5em}

Notice in this case $G/K^\pm\simeq S^4$. The cohomology groups of $M$ easily follow from the long exact sequence of the pair $(M,B_-)$.

\begin{conclusion*}
The manifold $M$ of type $N^7_I$ has the same homology groups as $S^4\times S^3$
\end{conclusion*}


%
%
%
%
%
%

\end{document}